\documentclass{scrartcl}

\usepackage[english]{babel}
\usepackage[utf8]{inputenc}
\usepackage[T1]{fontenc}
\usepackage{lmodern}

\usepackage{multirow}
\usepackage{rotating}
\usepackage{booktabs}
\usepackage{tabularx}
\newcolumntype{C}{>{\centering\arraybackslash}X}
\newcolumntype{R}{>{\raggedleft\arraybackslash}X}
\newcolumntype{L}{>{\raggedright\arraybackslash}X}

\usepackage{algorithm}
\usepackage{algorithmic}

\usepackage[usenames,dvipsnames,table]{xcolor}
\definecolor{rgbcyan}{HTML}{00ADEF}
\definecolor{rgbmagenta}{HTML}{EC008C}

\usepackage{tikz}
\usetikzlibrary{mindmap,backgrounds}
\usetikzlibrary{shapes,arrows}
\usetikzlibrary{matrix}
\usetikzlibrary{positioning}
\usetikzlibrary{calc}
\tikzstyle{line} = [draw, thick]
\usepackage[font=footnotesize,skip=4pt]{subcaption}

\graphicspath{{../}}

\usepackage{pifont}
\usepackage{amsmath}
\usepackage{amssymb}
\usepackage{amsthm}

\renewcommand{\vec}[1]{\boldsymbol{#1}}

\usepackage[breaklinks]{hyperref}
\usepackage{breakurl}

\usepackage{siunitx}

\title{Towards a Systematic Development \\Process of Optimization Methods}

\author{Simon Wessing\\Chair of Algorithm Engineering\\Computer Science Department\\Technische Universität Dortmund, Germany\\\texttt{simon.wessing@tu-dortmund.de}}

\date{}
\begin{document}
 
\maketitle

\begin{abstract}
The ultimate goal of all optimization methods is to solve real-world problems.
For a successful project execution, knowledge about optimization and the application has to be pooled.
As it is too inefficient to highly train one person in both fields, a team of experts is usually put together.
Unfortunately, communication errors must be expected when several people collaborate.
In this work, we deal with the avoidance and the repair of these communication errors.
The tools proposed in this regard are, among others, the algorithm engineering cycle,  checklists for structuring communication, and knowledge databases.
The discussion is enriched with examples from continuous optimization, but most tools have a much wider applicability, even beyond optimization.
\end{abstract}

\section{Introduction}

People working in algorithmics have long ago noticed the additional troubles caused by real-world problems compared to purely theoretical work~\cite{Weihe2001}. 
As a consequence, the field of \emph{algorithm engineering} was founded \cite{Mueller2010,Sanders2009}, to focus the attention of research in algorithmics more towards results with practical utility. 
The goal of the initiative was to 
establish a feedback loop in the algorithm development process, with experimental validation of theoretic results as the key innovation~\cite{Sanders2009}.
While this approach was accepted quickly and achieved considerable success~\cite{Chimani2010}, a ``blind spot'' still seems to have remained.

This blind spot is our ineptitude in dealing with troubles related to other professions, especially the social sciences.
For example, Weihe speaks of ``unavoidable communication problems'' between practitioners and algorithm developers~\cite{Weihe2001}.
The question arises if \emph{unavoidable} is really the right adjective.
In the author's opinion, the reason why so few progress has been made in this regard may be partly ignorance, because many of the difficulties probably are not considered scientifically interesting in the exact sciences, but more severe seems to be our inability to grasp what the actual problem with our communication is.

Especially in black-box optimization, where the author has personal experience, communication between the client and the optimization expert regarding the problem properties is often minimal. 
This seems to be due to the common advertising of black-box optimization algorithms as being able to solve a problem without any other information than the return values of the objective function. 
While this may be true in principle, it must always be emphasized that ideally the optimizer is already consulted for the problem formulation, to try to turn it into a ``grey box''.
As a rule of thumb, the more problem information is available, the better can the problem formulation and the optimization heuristic be custom tailored to each other and the easier gets the optimization.
Otherwise, clients may overestimate the abilities of optimization heuristics considerably.

In current practice, the optimizer often begins the development with little more information than the data type, dimension, and bounds of the search space and the knowledge that the problem is likely multimodal because local search procedure $X$ has been tried and yields unsatisfactory results.
Sometimes there even exists uncertainty about the ``natural'' data type of the search space. 
For example, software for embedded systems may be designed to avoid floating point arithmetic by using integers of an appropriately small measuring unit instead. 
So, is it advisable to apply a discrete method or a real-valued one?

Often project goals are missed, but not because the problem is overly difficult, but because the existing knowledge is not communicated properly. 
Thus, certain problem properties will be overlooked and consequently the optimization methods are applied wrongly. 
This ineptitude does not necessarily have to lead to a complete failure of the project, as it can be usually detected by a thorough analysis of the results.
However, repairing the situation is typically expensive, because the potentially laborious development process has to be reiterated, incorporating the newly gained information. 
For example, a ship propulsion optimization problem~\cite{Naujoks2007}, which was initially treated as a black box, turned out to have many more optima than expected~\cite{Rudolph2009}.
A satisfactory solution could finally be found with considerable effort~\cite{Rudolph2009}, but a conclusion was that it would have been advisable to spend more time in interdisciplinary teams to explore the problem and maybe find a less demanding formulation.

In his popular-science book ``The Checklist Manifesto'', Gawande~\cite{Gawande2011} proposes checklists as a useful tool for structuring routine tasks and teamwork in many professional fields. 
He argues that a major challenge of complex tasks is ``making sure we apply the knowledge we have consistently and correctly''~\cite[p.~10]{Gawande2011}.
His recommendations are based on his own experiences with the development of the world health organization's (WHO) safe surgery checklist~\cite{Haynes2009,Weiser2010}, which contains a small number of well-proven safety checks advisable for every operation.

It turns out that there is even an older example of checklists in a field very close to algorithm engineering.
This is the work of Coleman and Montgomery~\cite{Coleman1993} on planning of designed experiments. 
They propose \emph{pre-design experiment guide sheets}, which are essentially checklists in a \emph{read-do} form (meaning an action is taken after reading the corresponding item, as opposed to carrying out all tasks and confirming everything afterwards in a \emph{do-confirm} style).
They write that ``[t]he guide sheets are designed to be discussed and filled out by a multidisciplinary \emph{experimentation team} consisting of engineers, scientists, technicians/operators, managers, and process experts. 
These sheets are particularly appropriate for complex experiments and for people with limited experience in designing experiments.''

Based on these anecdotes, we identify three key topics for quality improvement in optimization:
\begin{enumerate}
\item systematic collection of existing expert knowledge,
\item better structuring of communication,
\item adopting an iterative improvement process.
\end{enumerate}
The core contribution of this work is to present these three topics as an integrated approach, with a focus on continuous optimization.
As a first step, we define a general algorithm development cycle in Sect.~\ref{sec:ae}.
The remaining sections then deal with measures to be taken to improve individual parts of this cycle.

\section{A Simplified Algorithm Engineering Cycle}
\label{sec:ae}

In the 1990s, a growing divide between theoretical and practical results was observed in algorithmics~\cite{Hooker1994}. 
As root causes, ignorance of characteristics of modern computer hardware and of real-world problem instances were identified~\cite{Chimani2010}.
The solution to these problems was the (re-)introduction of the scientific method for experimental research~\cite{Sanders2009}, establishing the field of algorithm engineering. 
Surveys on the topic can be found in~\cite{Chimani2010,Mueller2010,Sanders2009}. 
As all these texts have an intended audience of algorithm theorists, the carrying out of experiments is mostly presented as an extension to theoretical analysis. 
Additionally, there is a special emphasis on aspects related to (low-level) programming in all descriptions of algorithm engineering, as, e.\,g., consideration of memory hierarchies or software reuse by building libraries of algorithms and problems.

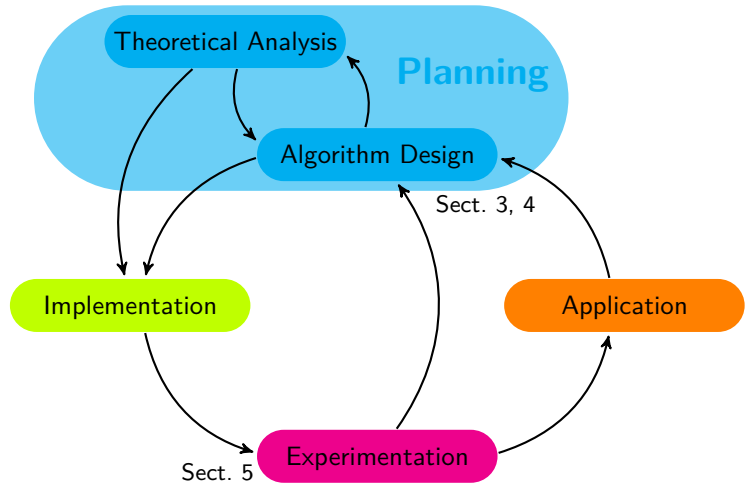
\begin{figure}[t]
\centering
\small
~~\begin{tikzpicture}
[text height=1.75ex,text depth=.25ex,->,>=stealth',shorten >=1pt,auto,font=\sffamily,
every node/.style={rectangle, fill=orange, text centered, rounded corners=1em, minimum height=2em, minimum width=9em, node distance=6em},
general/.style={fill=blue!20,thick},
detail/.style={thick},
onlytext/.style={thick,fill=none}]
\node[general,fill=cyan!50,minimum height=7em, minimum width=20em,rounded corners=3.5em] (plan) at (-1,2.75) {};

\node[detail,fill=cyan] (design) at (0,2) {Algorithm Design};
\node[detail,fill=cyan] (theory) at ($ (0,2) + (-2,1.5) $) {Theoretical Analysis}; 
\node[general,fill=magenta] (study) at ( 0,-2) {Experimentation};
\node[general,fill=orange] (act) at ( 3.25,0) {Application};
\node[general,fill=lime] (do) at (-3.25,0) {Implementation};

\node[onlytext] (planningsects) at (1.425,1.35) {\footnotesize Sect.~\ref{sec:determining_optproblem_checklist}, \ref{sec:expert_knowledge_optim}}; 
\node[onlytext] (expsec) at (-2.1,-2.2) {\footnotesize Sect.~\ref{sec:experimentation}}; 

\path[line] (do) edge [bend right, shorten >=1pt] (study)
(study)		edge [bend right=35] (design)
			edge [bend right] (act)
(design) edge [bend right, shorten >=0pt] (theory)
		 edge [bend right] (do)
(theory) edge [bend right, shorten >=0pt] (design)
		 edge [bend right] (do)
(act) edge [bend right] (design);
\node[general, fill=none, color=cyan] at (1.25,3.05) {\Large\textbf{Planning}};
\end{tikzpicture}
\caption{A simplified algorithm engineering cycle based on the one by Müller-Hannemann and Schirra~\cite{Mueller2010}.}
\label{fig:ae_cycle}
\end{figure}

Unfortunately, a rigorous analysis is seldom possible for inexact approaches such as heuristic optimization methods~\cite{Bartz-Beielstein2010}.
Müller-Hannemann and Schirra~\cite{Mueller2010} are the only ones who leave theoretical analysis out of the core development cycle and let it have its own separate cycle with algorithm design. 
(The two together correspond to the classical ``pen-and-paper'' approach to algorithmics~\cite{Sanders2009}.)
If we take~\cite{Mueller2010} as a basis and omit the extra topics mentioned above, we are left with the simplified development cycle in Fig.~\ref{fig:ae_cycle}. 
This cycle is also suited for purely experimental work and thus should be general enough to serve as a guideline for any algorithm development.
It suggests that theoretical analysis and real-world application are (at least in some iterations) optional to the development process.

While a feedback loop as in Fig.~\ref{fig:ae_cycle} should be standard practice in algorithm development nowadays, our goal is to reduce the number of costly cycles in it.
A first step towards this goal is a better up-front communication between the customer and the developer in the planning phase.
This is attempted in Sect.~\ref{sec:determining_optproblem_checklist}, by following a checklist in the same fashion as for designed experiments~\cite{Coleman1993}.
When the problem definition has been chosen, the algorithm design phase begins. 
Here, the developer has to take the information gathered in the meeting(s) and find a good optimization approach. 
At this point, it would be advantageous to have a lot of experience in optimization. 
However, in Sect.~\ref{sec:expert_knowledge_optim} we discuss strategies how to make life easier for us.
This encompasses collecting basic recommendations for less experienced people, making existing knowledge better searchable, and automating the development as much as possible.
Section~\ref{sec:experimentation} recapitulates some guidelines for good experimentation before everything is summarized in Sect.~\ref{sec:conclusion}.

\section{A Checklist for Determining the Optimization Problem}
\label{sec:determining_optproblem_checklist}

The attractivity of checklists lies in their ability to serve two purposes at once: guiding communication and storage of expert knowledge.
They have initially been developed in the aviation sector to structure the execution of routine tasks (see, e.\,g., Degani and Wiener~\cite{Degani1993,Degani1997}).
By now, the awareness of the virtues of checklists has also spread to other disciplines.
For example, there has recently been an initiative to introduce checklists into medicine in an attempt to make surgeries safer~\cite{Haynes2009,Weiser2010}. 
Earlier work of similar character could greatly reduce infections from central venous lines~\cite{Pronovost2006}. 
The benefit of checklists has also been recognized in statistics, where they are used for planning of designed experiments \cite{Coleman1993}. 

While the term ``checklist'' may for some people have a negative connotation of being rigid and restrictive, it should be noted that for good checklists this is not necessarily so. 
Instead of dictating the approach to solving a problem, they can also be used to specify \emph{communication tasks}~\cite[p.~65]{Gawande2011}. 
This means that the checklist simply requires the involved parties to discuss certain topics and to choose a solution that is deemed appropriate.

While checklists focus on executing tasks, questionnaires in comparison focus on data collection.
Thus, questionnaires may lay out several options, while checklists are rather phrased to define tasks and should only allow affirmative answers~\cite{Degani1993}.
If, however, the task is to answer a question, then the two forms can be obviously converted into one another.
For example, Doneit et al.~\cite{Doneit2015a} developed a questionnaire with a similar problem as ours in mind, namely to obtain information from experts about function approximation~\cite{Doneit2015}.
This form of presentation might be more appropriate when experts are to be interviewed individually, while checklists could be used more generally as an orientation guide for teams.

\begin{figure}[t!]
\hrule\smallskip
\small
\begin{minipage}{\textwidth}
\begin{enumerate}
\item Introduce participants, including their reason for attendance or, if already specified, their role in this project.
\item Define optimization goal (e.\,g., find feasible/robust/best solution, detect several local optima, approximate level set, approximate Pareto-set or -front)
\item State relevant background on objective function(s) and decision variables: (a) theoretical relationships; (b) expert knowledge/experience; (c) previous attempts/existing data. Where does this project fit into the study of the process or system?
\item List each objective function and answer the following questions:
	\begin{itemize}
	\item Can it be decomposed into terms of different meaning? Is it additively separable?
		\begin{itemize}
		\item If yes, treat each part separately in this discussion.
		\end{itemize}
	\item Is the analytic form/gradient information available?
	\item Is it linear/quadratic/convex/multimodal/unknown?
	\item If multimodal, does it have a global structure (i.\,e., funnel property, symmetries)?
	\item Is it deterministic? 
	\item What are its domain and image/range?
	\item What is its (distribution or range of) run time and/or cost?
	\end{itemize}
\item List each decision variable and answer the following questions:
	\begin{itemize}
	\item What is its domain (data type, lower/upper boundaries)?
	\item What is its current default value?
	\item What is the expected influence on each objective function?
	\item Should a nonlinear transformation (log, sqrt, etc.) be applied?
	\end{itemize}
\item List each side constraint and answer the following questions:
	\begin{itemize}
	\item Is it known (K) or hidden (H)?
	\item Is it an a priori (A) constraint or simulation-based (S)? Related to this: What is its (distribution or range of) run time and/or cost?
	\item Is it relaxable (R) or unrelaxable (U)?
	\item Is it quantifiable (Q) or nonquantifiable (N)?
	\end{itemize} 
\item If more than one objective function, list possibly conflicting pairs.
\item Decide on a problem formulation, potentially limited to a subset of the components identified in steps 4 to 6.
\item Choose a cost model and the (maximum) budget to spend on the problem.
\item Allocate responsibilities to people.
\end{enumerate}
\end{minipage}
\smallskip\hrule
\caption{A checklist for determining optimization problems.}
\label{fig:opt_checklist}
\end{figure}

Figure~\ref{fig:opt_checklist} contains a checklist for determining an optimization problem.
To save space, it is only presented in condensed form, without check boxes and space for the answers.
It should be used in the inaugural project meeting (and potentially in following meetings that mark the beginning of a new iteration in the development cycle, see Fig.~\ref{fig:ae_cycle}) to obtain a common understanding of the problem and to ensure that all the expert knowledge is transferred.
This hopefully would help to prevent avoidable mistakes, which might otherwise be made due to misunderstandings or omissions. 
The checklist's style closely resembles that of~\cite{Coleman1993} and likewise it should be completed by the client and the optimizer together. 
By working through the questions, the answers should be written down, either on paper or digitally.
Note that it may not be possible to answer all questions completely.

We will now explain each item in detail.

\paragraph{Introductory Part}
The first item requires all participants of the meeting to introduce themselves, in aspiration to enable good teamwork. 
(This idea is taken from the WHO checklist~\cite{Haynes2009}.) 
The goal of the second item is to ask the customers what their actual requirement is for the outcome of the optimization. 
This should especially reveal if an a-priori approach is sought, where the optimization algorithm returns a single solution, an a-posteriori method, where a set of solutions is prepared for human inspection or some other post-processing~\cite{Miettinen2008}, or even an interactive method with a human in the loop. 
The intention of item three is to inquire about general expert knowledge that may be available. 
Especially, it should be figured out which degree of maturity the studied system and its current solution has. 
This may indicate how much effort is required for further improvement, according to the Pareto principle.
Furthermore, the context of the optimization problem should be explained.

After the initial questions have been answered, there come three sections regarding objective functions, decision variables, and side constraints. 
The order in which these sections are processed does not seem too important and for some problem formulations they may even be intertwined.
Note that the sets of the three entities collected here should not be pruned prematurely, but should contain all \emph{potentially} useful candidates for objectives, variables, and constraints.

\paragraph{Objectives}
Regarding the objective functions, several properties must be checked. 
First of all, it is inquired if the function is composed of more elementary functions, e.\,g., if $f(x,y)$ can be represented as $g(x,y) + h(x,y)$.
This is an important issue, because the sub-functions are expected to have a different meaning and thus could be conflicting.
A special case of composition is an additively separable function $f(x, y) = g(x) + h(y)$, which could be divided into completely independent sub-problems.
In any case, the identified components of the functions should be treated separately in following questions.

The most important properties are retrieved in the following.
A very important one is if the analytic form of the function is known, because this has severe consequences on the preferred optimization approach. 
If the analytic form is known and it is differentiable, then a gradient descent method, quasi-Newton method, or maybe even an analytic solution could be possible.
To this end, the type of the function should be narrowed down further.
If the type is unknown or the function is multimodal, further (possibly estimated) properties become relevant.
For example, important discrete and real-valued multimodal problems exhibit a certain global structure in which local optima are correlated in quality and position. 
Examples are the traveling salesperson problem~\cite{Lourenco2010} or molecular conformation problems~\cite{Wales1999}.
It is well known that if a problem has such a so-called funnel structure, it is advisable to use certain specialized approaches such as iterated local search (discrete case)~\cite{Lourenco2010}, basin hopping~\cite{Wales1999}, or IPOP-CMAES (real-valued case)~\cite{Auger2005}. 
However, in reality the property is not binary and so far there exists no approach to quantify the degree of ``funnelishness'' of a problem. 
So, the current practice is mainly to guess the property's value from the observed performance of the corresponding optimization algorithms~\cite{Addis2011}.

Another relevant feature would be random noise in the evaluations, which can stem from various sources~\cite{Beyer2007}. 
Which decision variables form the domain of the function is of course mandatory information and also to know lower or upper bounds for the function values can be useful.
These properties naturally have to be given for each individual function, because, e.\,g., the domains of different objective functions naturally do not have to be identical~\cite[p.~125]{DagstuhlReport15031}. 
Finally, also the run time and other costs of the objective function should be estimated. 
The costs may be considered either constant or variable (e.\,g., dependent on some solution properties).
If possible, also further theoretical and/or empirical properties should be listed. 
To this end, it may also be useful to go through the questionnaire of Doneit et al.~\cite{Doneit2015a}.

\paragraph{Decision Variables}
For decision variables, data type and boundaries should be given. 
If no bounds are known, this may indicate a client's inexperience with the problem.
If an implemented solution exists, it would be interesting what the currently used values are in the decision variables. 
Furthermore, the expert should give an estimate how influential a decision variable is on each objective function. 
This would affect the decision which variables to select for the final problem formulation and possibly the chosen optimization approach, too. 
Sometimes, if a decision variable spans several orders of magnitude, it can be very beneficial to apply a nonlinear transformation~\cite{Preuss2015,Wagner2012}, e.\,g., the logarithm.

\paragraph{Constraints}
Side constraints are somewhat similar to objective functions. 
Here, too, it is most important how much information we can get, i.\,e., is the error condition specified, do we have a mathematical formulation of the constraint, can we compute the amount of constraint violation? 
The questions in this section constitute the QRAK classification scheme of Le Digabel and Wild~\cite{LeDigabel2015}.
In this scheme, each constraint is described by a four-letter code, which indicates the difficulty introduced by this constraint and can be linked to directions how to best treat it. 
For example, a hidden constraint must be necessarily also classified as nonquantifiable, unrelaxable, and simulation-based (NUSH). 
It is the worst kind of constraint, because it essentially means dealing with faulty software. 
Known constraints can be categorized into eight classes with the remaining questions.
The a priori group corresponds to relatively ``cheap'' constraints, where usually standard approaches exist for certain subclasses (e.\,g., linear, nonlinear). 
Simulation-based constraints, on the other hand, require potentially expensive calls to simulation software. 
A constraint is called relaxable if the objective function can still be reliably evaluated even when the constraint is violated, and unrelaxable otherwise.
Finally, quantifiability means that the degree of feasibility and/or violation can be quantified.
For further explanation and examples, we refer to the original paper~\cite{LeDigabel2015}.

\paragraph{Concluding Tasks}
To complete the checklist, the team has to decide on a final problem formulation. 
To do this, it seems advisable to discuss which objective functions are conflicting and which are not, so the number of objectives in the final problem can be kept low. 
Some functions also may be switched from objective to constraint or vice versa.
(Note that it is generally not unusual to treat the optimization problem under a different paradigm than the most obvious one.
For example, single-objective problems can be treated by multiobjective methods~\cite{Preuss2015,Segura2013a}, multiobjective problems can be treated as single-objective by scalarization~\cite{Klamroth2013}, and so on.)
Then it has to be decided which decision variables make up the search space and which are held constant. 
The cost model, e.\,g., time or some bottleneck operation, has to be fixed and the budget must be determined.
Finally, the contributions and responsibilities of the people involved in the project are set out in writing.

\subsection{Remarks on the Design and Validation of Checklists}
As situations in medicine and aviation are often very time-critical, the prevailing checklist style there is brief and concise.
As we are not subject to such restrictions, also more extensive variants are possible. 
Thinking this further, they could also be extended to decision trees, as Gunter proposes in~\cite{Gunter1993}. 
However, the format should be such that everything can be worked through in one meeting.

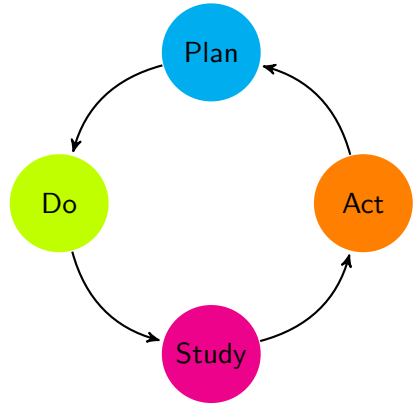
\begin{figure}[t]
\centering
\begin{tikzpicture}
[text height=1.75ex,text depth=.25ex,->,>=stealth',shorten >=1pt,auto,font=\sffamily,
general/.style={circle,fill=blue!20,thick,
inner sep=0pt,minimum size=13mm},
detail/.style={rectangle,draw=black!50,fill=black!20,thick,
inner sep=0pt,minimum size=4mm}]
\node[general,fill=cyan] (plan) at ( 0,2) {Plan};
\node[general,fill=magenta] (study) at ( 0,-2) {Study};
\node[general,fill=orange] (act) at ( 2,0) {Act};
\node[general,fill=lime] (do) at (-2,0) {Do};

\path[line] (plan) edge [bend right] (do)
(do) edge [bend right] (study)
(study) edge [bend right] (act)
(act) edge [bend right] (plan);
\end{tikzpicture}
\caption{The plan-do-study-act cycle. Note the similarity to Fig.~\ref{fig:ae_cycle}.}
\label{fig:pdsa}
\end{figure}

Weiser et al.~\cite{Weiser2010a} report on their approach of developing the WHO surgical safety checklist. 
After an initial collection of potentially relevant topics and converting them into checklist items, points were sought in the surgical workflow where interruption by reading a checklist seemed fitting.
The initial checklist was then iteratively refined by small trials, using a process oriented to the plan-do-study-act (PDSA) paradigm~\cite{Moen2010}.
The PDSA cycle is based on the scientific method~\cite{Moen2010} and was popularized as a quality improvement process in the 20th century by Deming~\cite[p.~88]{Deming1988}. 
Figure~\ref{fig:pdsa} illustrates its four stages: 1)~plan: define the objectives and changes to achieve them, 2) do: implement the changes, 3) study: analyze the results, 4) act: decide whether to accept or reject the changes; if accepted, carry out changes on a large scale, making them the new standard.
Due to the shared origin in the scientific method, PDSA and the algorithm engineering cycle are quite similar.

Weiser et al.\ finally validated the checklist by conducting a large empirical study in several hospitals around the world, obtaining positive results~\cite{Weiser2010}.
Also checklists in aviation are developed by iteratively refining them in flight simulators.
While the other checklist and questionnaire examples~\cite{Coleman1993,Doneit2015a} may have been designed by domain experts, they have not been experimentally validated, to the authors best knowledge.
Unfortunately it also seems similarly laborious to validate checklists in optimization as it is difficult to recruit a sufficient number of qualified test persons and find a large number of applications.
Additionally, checklists for the structuring of communication probably depend largely on culture and personal preference, as it is often mentioned that they may be individually adapted to it~\cite{Coleman1993,Weiser2010a}.
In summary, it does not seem too severe if communication-oriented checklists are not immediately validated in a controlled environment. 
The proposal in Fig.~\ref{fig:opt_checklist} shall serve as a stimulus to start a discussion and exploration of this topic.

Inspiration on possible approaches for future validation may be taken from software engineering, where the empirical comparison and assessment of programming languages and tools is also very challenging, but slowly making progress~\cite{Oram2010}.

\section{Expert Knowledge in Optimization}
\label{sec:expert_knowledge_optim}

In this section we discuss possible ways to obtain, store, and distribute knowledge in optimization, to reduce the dependency on the developer's capability in the algorithm design step (see Fig.~\ref{fig:ae_cycle}).

\subsection{Basic Recommendations}
\label{sec:basic_recommendations}

For people with limited experience, with a tight schedule, or without academic requirements, some general advice in abridged form is certainly helpful for quick results.
This advice may consist of simple rules of thumbs, which proved beneficial in most situations.
For example, the safety checks in the WHO checklist~\cite{Haynes2009} are of this character.
Also in optimization, some rules of thumb have almost universal validity, an example of which is given in Sect.~\ref{sec:example_search_space_scaling}.
Checklists or decision trees again seem to be a good form of presentation.
They could be used as model for choosing a general optimization approach for a given problem, as designing a new algorithm from scratch should be usually avoided if possible.
A first attempt into this direction is taken by Mittelmann~\cite{Mittelmann}.
Also the MathWorks documentation webpages for the Optimization Toolbox~\cite{MathWorks} and the Global Optimization Toolbox~\cite{MathWorksa} contain decision tables providing advice which algorithm to use, depending on the problem class.

Unfortunately, further advice when to use which method is widely scattered in the literature, and only few papers explicitly discuss strengths \emph{and weaknesses} of methods.
This has of course partly political reasons, but often it is also because the experiments carried out are not large enough to provide such information (more on this topic in Sect.~\ref{sec:experimentation}).
Some positive examples for papers containing general advice in continuous optimization are~\cite{Forrester2009,Jones2001,Torczon1997,TorAliVii99}.

\subsubsection{An Example from Continuous Optimization}
\label{sec:example_search_space_scaling}

Some pitfalls of optimization problems are frequently neglected in the scientific literature, because the workarounds are simple.
An example for this are implicit assumptions about the size of the search space often built into optimization methods in continuous optimization. 
Some algorithms handle this topic explicitly by requiring the specification of bound constraints as part of the problem description. 
Examples for this group are model-based optimization algorithms~\cite{Jones2001}. 
They use this information to draw a uniform initial sample and to bound the space in which to search for promising points for improvement of the model.
Local search algorithms, e.\,g., direct search methods, do not necessarily have this requirement, because they only use a single starting point. 
However, people should recognize that their initial step size contains an implicit assumption about the size of the attraction basin it is started in. 
(We would rather not choose the step size larger than the expected distance to the optimum we are interested in.)
If this assumption is wrong, i.\,e., especially if the step size is too small, the obtained performance may be very bad even on unimodal functions. 
Here are two examples:
\begin{itemize}
\item The simplex search by Nelder and Mead~\cite{Nelder1965} is a direct search method. 
For its search in $n$ dimensions, the algorithm uses a simplex consisting of $n + 1$ points. 
Most of the common implementations today provide no interface for specifying a region of interest, an initial step size, or an initial simplex.
Instead, they construe the initial simplex from the starting point by using questionable heuristics.
In Matlab's fminsearch and SciPy, 
the points $\vec{x}_i$, $i = 2, \dots, \mbox{n + 1}$ of the simplex are obtained from the starting point $\vec{x}_1$ by the following rule:
\[
\vec{x}_i = \begin{cases}
(x_{1,1}, \dots, x_{1,i-2}, 1.05x_{1,i-1}, x_{1,i}, \dots, x_{1,n})^\top  & \text{if } x_{1,i-1} \neq 0\,,\\
(x_{1,1}, \dots, x_{1,i-2}, 0.00025, x_{1,i}, \dots, x_{1,n})^\top & \text{else.}\\
\end{cases}
\]
This approach is purportedly due to an ``L. Pfeffer at Stanford'', according to \cite{Fan2002}.
As a consequence of this approach, a very small simplex is generated if the starting point is close to the origin. 
Even if only one coordinate is close to zero, the simplex is more or less degenerated, which is also bad. 
The implementation in R's optim method is slightly better, using
\[
\vec{x}_i = 
(x_{1,1}, \dots, x_{1,i-2}, x_{1,i-1} + 0.1\cdot\max\{|x_{1,1}|, \dots, |x_{1,n}|\}, x_{1,i}, \dots, x_{1,n})^\top 
\]
as construction rule.
This rule, which was adopted from \cite[pp.~168-178]{Nash1990}, is only problematic if \emph{all} coordinates of the starting point are close to zero. 
However, note that there is in general no reason to assume such a special role for the origin in the search space.
\item The recent survey of optimization algorithms by Rios and Sahinidis~\cite{Rios2013} uses an unusually large search space of $[-10^4, 10^4]^n$ for the benchmarking. 
This search space normally does not match the assumptions of algorithms' default parameters. 
Some of the algorithms in the competition were specifically adapted to this search space, and one of these performed best in the comparison. 
However, this is merely evidence for a good parameter tuning and not so much for the algorithm's superiority over the other contestants.
\end{itemize}
Normalizing the search space (or the region the starting points are drawn from) to $[0, 1]^n$ would eliminate the latter problem.
There do not seem to be any drawbacks associated with this action, so we propose to always do it in continuous optimization. 
Other authors have given this recommendation before, but probably did not insist strongly enough~\cite{Addis2011,Battermann2002}. 
Additionally, developers of unconstrained optimization algorithms should make sure that their algorithm performs best in the unit hypercube, e.\,g., by choosing an appropriate default step size. 
This is already often done, because floating point operations after all yield the highest precision close to zero, but it should be communicated more explicitly. 
Additionally, the user should have the opportunity to modify the initialization, to handle special cases.

\subsection{Networked Science}

In Sect.~\ref{sec:basic_recommendations}, several papers were cited for the useful advice they were giving.
However, the corresponding evidence is not (completely) included in these papers.
This is of course unsatisfactory from a scientific point of view, but is a general problem of the relatively short articles forming our main communication channel.
Another difficulty of this format is to aggregate results from different persons and to update existing knowledge with new information~\cite{Crepinsek2014}.

The question arises how to store and communicate data that is too bulky or not scientifically interesting enough to be included in publications.
To the author's best knowledge, large-scale comparisons in optimization are currently mostly carried out in competitions and special sessions at conferences.
In machine learning, the topic seems already more advanced. 
People there have begun to tackle the problem by storing raw data in experiment databases~\cite{Blockeel2007,Vanschoren2008} and subsequently building a website around these databases~\cite{Vanschoren2014}. 
The website is now organized around the basic concepts of machine learning experiments:
a small number of pre-defined machine learning \emph{tasks} (e.\,g., classification, regression) can be associated with the input data sets, \emph{workflows} describe algorithms to solve certain tasks, and \emph{runs} contain the results of executions of such workflows. 
Anybody can register at the website and upload their own work, while the website presents everything in an attractive way to the general public.
The whole approach nowadays runs under the catch phrase \emph{networked science}~\cite{Vanschoren2014}.
It seems equally useful for optimization, but setting up such a project is of course extremely labor-intensive.

\subsection{Automated Algorithm Design}
\label{sec:coseal}

The practical suitability of any algorithm depends on problem features.
Using inductive reasoning, one could optimistically assume that observed properties of a sample of instances from a certain problem class are also given for other instances of the problem.
This reasoning is inherently uncertain, but often successful in practice.
For example, in algorithm engineering there exist many reports of applications with special distributions of problem instances, which only make up a small fraction of the space of theoretically possible instances for the corresponding problem class.
This expert knowledge has been successfully exploited in ``manual'' algorithm development~\cite{Chimani2010}, by specially adapting an algorithm to the sample of instances.
From a point of view of machine learning, the adapted algorithm constitutes a constant prediction model, unconditionally choosing the same algorithm for all future problem instances.  

It is possible to reduce human involvement in the design process, by treating the search for the best algorithm as an optimization problem.
This problem is generally known under the terms \emph{algorithm selection}~\cite{Bischl2016} or \emph{algorithm configuration}~\cite{Hutter2014}, depending on whether the emphasis is on the choice between different algorithms or between different parameter settings of the same algorithm. 
(However, note that this distinction is rather philosophical.)
Unfortunately, such problems are typically expensive and noisy. 
Algorithm configuration often also deals with mixed variables, which is why algorithm selection is normally regarded as the easier of the two.

Usually, a model for the configuration/selection is built during an offline training phase. 
The simplest case of a constant model was already mentioned above.
If we want a conditional model, appropriate problem features must be available. 
The systematic search for such features is nowadays a research topic of its own~\cite{Mersmann2013}.
In continuous optimization, it is called \emph{exploratory landscape analysis} (ELA)~\cite{MerBisTra11,MerPreTra11}.
Based on these features, we may then predict an algorithm's performance on a problem instance, and consequently choose the best algorithm from a portfolio.
The approach is useful when an algorithm can be identified, being so much better than the default algorithm that it can make up for the additional cost of ELA. 
Thus, calculating the features should be as cheap as possible.

Apart from a conventional offline training approach, algorithm selection may also be treated as a multi-armed bandit problem~\cite{Gagliolo2010}.
In this case, no features are required \linebreak at all.

\section{Experimentation}
\label{sec:experimentation}

In comparison to the planning of optimization projects, the experimental analysis of optimization algorithms is nowadays a relatively well-researched topic. 
Twenty years ago, this was not so.
Experimental analyses were less accepted than today~\cite{Hooker1994}, and existing work was criticized as unsuitable to yield insight into the working principles of algorithms~\cite{Hooker1995}. 
As a root cause, Hooker identified the competitive nature of most experiments, which in the long run leads to overfitting to the used benchmark sets~\cite{Hooker1995}.
He contrasts this approach with \emph{scientific testing}, which requires to incorporate more factors into the experimental design and analysis, to actually be able to explain the made observations. 

\begin{figure}[t]
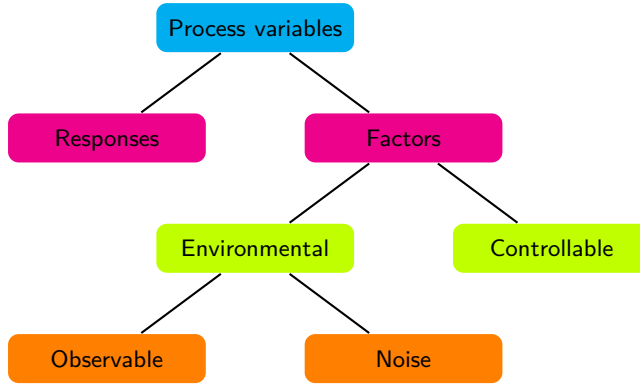
%
\centering
\footnotesize
\sffamily
\tikz [thick, every node/.style={rectangle, fill=orange, text centered, rounded corners, minimum height=2em, minimum width=8em,text height=1.75ex,text depth=.25ex}, level 1/.style={sibling distance=12em},
level 3/.style={sibling distance=12em}, level distance=4.5em] 
\node [fill=cyan] {Process variables} 
	child { node [fill=magenta] {Responses} }
	child { node [fill=magenta] {Factors}
		child { node [fill=lime] {Environmental} 
			child { node {Observable} }
			child { node {Noise} }
		}
		child { node [fill=lime] {Controllable} }
	};
\caption{Different categories of process variables.}
\label{fig:process_variables}
\end{figure}

To understand what this means, we have to give a brief introduction to the terminology of designed experiments (also see \cite[pp.~3--16]{Castillo2007}).
A \emph{designed} experiment is an experiment in which the test locations are planned by the experimenter, and this set of locations is called the experimental design.
The test locations are points in the region of interest, which is in turn a subset of the possible input space of the process or system in consideration.
The input variables are called factors. 
Responses are outputs of the process.
In his proposal for scientific testing, Hooker essentially suggests to not only use one factor ``optimization algorithm'', but to identify modular components of algorithms.
Furthermore, it was requested several times to incorporate problem features as factors in the experiment~\mbox{\cite{Eiben2002,Hooker1994}}.

At this point, it seems appropriate to distinguish three different types of factors: controllable, observable, and non-observable/noise (in~\cite{Coleman1993} the latter are called nuisance factors).
This classification is shown in Fig.~\ref{fig:process_variables}.
It is not novel, but also not standard practice in optimization. 
Thus, many benchmarking experiments suffer from confusion about which information is available in the real-world.
The classification is justified as follows: 
naturally, the planning of test locations is only possible for controllable factors, e.\,g., the properties of optimization algorithms.
In real-world applications, some of the factors cannot be controlled or even observed.
This applies especially to most, but not all features of optimization problems. 
The ``observability'' of a factor is important, because it has implications for the factor's treatment in real-world applications. 
For example, a problem's number of decision variables $n$ is always known to us and therefore an optimization algorithm is not required to possess a good performance over all possible values of~$n$. 
Instead we can always choose the most appropriate algorithm for the current $n$ before we begin the optimization (given that we somehow possess this expert knowledge). 
On the other hand, tuning or otherwise configuring an algorithm for a specific problem instance and reporting the tuned results without accounting for the additional cost is unsound, if the problem instance is a black box in reality.

Luckily we can actually control most if not all of the factors (even pseudorandom numbers) in computational experiments using artificial problems.
This general observation is usually ascribed to Taguchi (see, e.\,g., \cite[pp.~223--234]{Castillo2007} or \cite{SacWelMit89}), who proposed to incorporate uncontrollable factors into the experimental design, to obtain results exhibiting robustness later in the real world. 
The general approach of the experimental evaluation should be to select the best configuration among the control factors, regarding all possible levels of noise factors, but depending on particular levels of the observable factors.

Unsurprisingly, we recommend to use Coleman and Montgomery's checklist~\cite{Coleman1993} for planning the experiment. 
It should especially help to classify the factors according to Fig.~\ref{fig:process_variables}.
Additionally, a standardized scheme should be used for reporting on the experiment. 
Barr et al.~\cite{Barr1995} give some general advice in this regard.
Preuss~\cite{Bartz-Beielstein2010,Preuss2007} more concretely proposes an organization into research question, pre-experimental planning, task, setup, results, observations, and discussion. This structure aims to facilitate the distinction of objective results and subjective comments, and to improve the reproducibility of results without access to the implementation.
Regarding the analysis of the obtained data, Montgomery emphasizes the importance of visualizations in his practical guidelines, especially for presenting the results to others \cite[p.~16]{Montgomery1997}. 
Trellis graphics are a viable method to investigate on the interactions of the various factors~\cite{Becker1996}.

\section{Summary}
\label{sec:conclusion}
In this work we suggested concrete measures to improve the development of optimization algorithms, especially with regard to real-world applications.
The basis of the proposal is a simplified algorithm engineering cycle (cf.~Fig.~\ref{fig:ae_cycle}), which serves as an iterative development process.
This cycle consists of planning, implementation, experimentation, and application. 
For the planning stage, a checklist was proposed to guide an interdisciplinary team in defining the optimization problem. 
The purpose of the checklist is not to prescribe solutions, but to structure communication.
Furthermore, the research community was encouraged to explore new ways of storing and exchanging results, apart from research papers.
The algorithm design, which was regarded as a part of the planning stage here, could be automated more, leading in principle to another optimization problem in the algorithm parameter space. 
Finally, we briefly dealt with experiments, recommending to classify factors according to their observability in the real-world, to eventually increase the information content of experimental results.
Also the reporting on experiments should be structured more.

The implementation of optimization algorithms was not discussed in this work, as the topic should be sufficiently covered by the algorithm engineering~\cite{Moerig2010} and general software development literature~\cite{Oram2010}. 
Also the application was not dealt with because it seems completely case-dependent.

\bibliographystyle{plain}
\bibliography{checklist_literature.bib}   

\begin{thebibliography}{10}

\bibitem{Addis2011}
Bernardetta Addis, Andrea Cassioli, Marco Locatelli, and Fabio Schoen.
\newblock A global optimization method for the design of space trajectories.
\newblock {\em Computational Optimization and Applications}, 48(3):635--652,
  2011.

\bibitem{Auger2005}
Anne Auger and Nikolaus Hansen.
\newblock A restart {CMA} evolution strategy with increasing population size.
\newblock In {\em {IEEE} Congress on Evolutionary Computation ({CEC})},
  volume~2, pages 1769--1776, 2005.

\bibitem{Barr1995}
Richard~S. Barr, Bruce~L. Golden, James~P. Kelly, Mauricio~G.C. Resende, and
  William~R. Stewart, Jr.
\newblock Designing and reporting on computational experiments with heuristic
  methods.
\newblock {\em Journal of Heuristics}, 1(1):9--32, 1995.

\bibitem{Bartz-Beielstein2010}
Thomas Bartz-Beielstein and Mike Preuss.
\newblock The future of experimental research.
\newblock In Thomas Bartz-Beielstein, Marco Chiarandini, Luís Paquete, and
  Mike Preuss, editors, {\em Experimental Methods for the Analysis of
  Optimization Algorithms}, pages 17--49. Springer, 2010.

\bibitem{Battermann2002}
Astrid Battermann, Joerg~M. Gablonsky, Alton Patrick, Carl~T. Kelley,
  Kathleen~R. Kavanagh, Todd Coffey, and Cass~T. Miller.
\newblock Solution of a groundwater control problem with implicit filtering.
\newblock {\em Optimization and Engineering}, 3(2):189--199, 2002.

\bibitem{Becker1996}
Richard~A. Becker, William~S. Cleveland, and Ming-Jen Shyu.
\newblock The visual design and control of trellis display.
\newblock {\em Journal of Computational and Graphical Statistics},
  5(2):123--155, 1996.

\bibitem{Beyer2007}
Hans-Georg Beyer and Bernhard Sendhoff.
\newblock Robust optimization – a comprehensive survey.
\newblock {\em Computer Methods in Applied Mechanics and Engineering},
  196(33-34):3190--3218, 2007.

\bibitem{Bischl2016}
Bernd Bischl, Pascal Kerschke, Lars Kotthoff, Marius Lindauer, Yuri Malitsky,
  Alexandre Fr{\'{e}}chette, Holger~H. Hoos, Frank Hutter, Kevin Leyton-Brown,
  Kevin Tierney, and Joaquin Vanschoren.
\newblock {ASlib}: A benchmark library for algorithm selection.
\newblock {\em Artificial Intelligence}, 237:41--58, 2016.

\bibitem{Blockeel2007}
Hendrik Blockeel and Joaquin Vanschoren.
\newblock Experiment databases: Towards an improved experimental methodology in
  machine learning.
\newblock In Joost~N. Kok, Jacek Koronacki, Ramon Lopez~de Mantaras, Stan
  Matwin, Dunja Mladenič, and Andrzej Skowron, editors, {\em Knowledge
  Discovery in Databases: PKDD 2007}, volume 4702 of {\em Lecture Notes in
  Computer Science}, pages 6--17. Springer, 2007.

\bibitem{Chimani2010}
Markus Chimani and Karsten Klein.
\newblock Algorithm engineering: Concepts and practice.
\newblock In Thomas Bartz-Beielstein, Marco Chiarandini, Luís Paquete, and
  Mike Preuss, editors, {\em Experimental Methods for the Analysis of
  Optimization Algorithms}, pages 131--158. Springer, 2010.

\bibitem{Coleman1993}
David~E. Coleman and Douglas~C. Montgomery.
\newblock A systematic approach to planning for a designed industrial
  experiment.
\newblock {\em Technometrics}, 35(1):1--12, 1993.

\bibitem{Degani1993}
Asaf Degani and Earl~L. Wiener.
\newblock Cockpit checklists: Concepts, design, and use.
\newblock {\em Human Factors}, 35(2):345--359, 1993.

\bibitem{Degani1997}
Asaf Degani and Earl~L. Wiener.
\newblock Procedures in complex systems: the airline cockpit.
\newblock {\em IEEE Transactions on Systems, Man and Cybernetics, Part A:
  Systems and Humans}, 27(3):302--312, 1997.

\bibitem{Castillo2007}
Enrique del Castillo.
\newblock {\em Process Optimization}, volume 105 of {\em International Series
  in Operations Research \& Management Science}.
\newblock Springer, 2007.

\bibitem{Deming1988}
W.~Edwards Deming.
\newblock {\em Out of the crisis}.
\newblock Cambridge University Press, 1986.

\bibitem{Doneit2015a}
Wolfgang Doneit, Ralf Mikut, Lutz Gröll, and Markus Reischl.
\newblock {Fragebogen zur Erfassung von Vorwissen in Funktionsapproximationen}.
\newblock Technical report, Institut für Angewandte Informatik, Karlsruher
  Institut für Technologie, 2015.
\newblock (in German) \url{https://dx.doi.org/10.13140/RG.2.1.3511.3446}.

\bibitem{Doneit2015}
Wolfgang Doneit, Ralf Mikut, Lutz Gröll, and Markus Reischl.
\newblock {Vorwissen in Funktionsapproximationen durch
  Support-Vektor-Regression bei schlechter Datenqualität}.
\newblock In Frank Hoffmann and Eyke Hüllermeier, editors, {\em Proceedings
  25. Workshop Computational Intelligence}, volume~54 of {\em Schriftenreihe
  des Instituts für Angewandte Informatik / Automatisierungstechnik}, pages
  163--181. KIT Scientific Publishing, 2015.
\newblock (in German).

\bibitem{Eiben2002}
Agoston~E. Eiben and Mark Jelasity.
\newblock A critical note on experimental research methodology in {EC}.
\newblock In {\em IEEE Congress on Evolutionary Computation (CEC)}, volume~1,
  pages 582--587, 2002.

\bibitem{Fan2002}
Ellen Fan.
\newblock Global optimization of {Lennard-Jones} atomic clusters.
\newblock Master's thesis, McMaster University, 2002.

\bibitem{Forrester2009}
Alexander~I.J. Forrester and Andy~J. Keane.
\newblock Recent advances in surrogate-based optimization.
\newblock {\em Progress in Aerospace Sciences}, 45(1-3):50--79, 2009.

\bibitem{Gagliolo2010}
Matteo Gagliolo and Catherine Legrand.
\newblock Algorithm survival analysis.
\newblock In Thomas Bartz-Beielstein, Marco Chiarandini, Luís Paquete, and
  Mike Preuss, editors, {\em Experimental Methods for the Analysis of
  Optimization Algorithms}, pages 161--184. Springer, 2010.

\bibitem{Gawande2011}
Atul~A. Gawande.
\newblock {\em The Checklist Manifesto}.
\newblock Profile Books, 2011.

\bibitem{DagstuhlReport15031}
Salvatore Greco, Kathrin Klamroth, Joshua~D. Knowles, and Günter Rudolph.
\newblock Understanding complexity in multiobjective optimization {(Dagstuhl
  Seminar 15031)}.
\newblock {\em Dagstuhl Reports}, 5(1):96--163, 2015.

\bibitem{Gunter1993}
Berton~H. Gunter.
\newblock [{A} systematic approach to planning for a designed industrial
  experiment]: Discussion.
\newblock {\em Technometrics}, 35(1):13--14, 1993.

\bibitem{Haynes2009}
Alex~B. Haynes, Thomas~G. Weiser, William~R. Berry, Stuart~R. Lipsitz,
  Abdel-Hadi~S. Breizat, E.~Patchen Dellinger, Teodoro Herbosa, Sudhir Joseph,
  Pascience~L. Kibatala, Marie Carmela~M. Lapitan, Alan~F. Merry, Krishna
  Moorthy, Richard~K. Reznick, Bryce Taylor, and Atul~A. Gawande.
\newblock A surgical safety checklist to reduce morbidity and mortality in a
  global population.
\newblock {\em New England Journal of Medicine}, 360(5):491--499, 2009.
\newblock PMID: 19144931.

\bibitem{Hooker1994}
John~N. Hooker.
\newblock Needed: An empirical science of algorithms.
\newblock {\em Operations Research}, 42(2):201--212, 1994.

\bibitem{Hooker1995}
John~N. Hooker.
\newblock Testing heuristics: We have it all wrong.
\newblock {\em Journal of Heuristics}, 1(1):33--42, 1995.

\bibitem{Hutter2014}
Frank Hutter, Manuel López-Ibáñez, Chris Fawcett, Marius Lindauer, Holger~H.
  Hoos, Kevin Leyton-Brown, and Thomas Stützle.
\newblock {AClib}: A benchmark library for algorithm configuration.
\newblock In Panos~M. Pardalos, Mauricio~G.C. Resende, Chrysafis Vogiatzis, and
  Jose~L. Walteros, editors, {\em Learning and Intelligent Optimization},
  volume 8426 of {\em Lecture Notes in Computer Science}, pages 36--40.
  Springer, 2014.

\bibitem{Jones2001}
Donald~R. Jones.
\newblock A taxonomy of global optimization methods based on response surfaces.
\newblock {\em Journal of Global Optimization}, 21(4):345--383, 2001.

\bibitem{Klamroth2013}
Kathrin Klamroth, Elisabeth Köbis, Anita Schöbel, and Christiane Tammer.
\newblock A unified approach for different concepts of robustness and
  stochastic programming via non-linear scalarizing functionals.
\newblock {\em Optimization}, 62(5):649--671, 2013.

\bibitem{LeDigabel2015}
Sebastien Le~Digabel and Stefan~M. Wild.
\newblock A taxonomy of constraints in simulation-based optimization.
\newblock Technical Report ANL/MCS-P5350-0515, Argonne National Laboratory,
  2015.
\newblock \url{http://www.mcs.anl.gov/papers/P5350-0515.pdf}.

\bibitem{Lourenco2010}
Helena~R. Lourenço, Olivier~C. Martin, and Thomas Stützle.
\newblock Iterated local search: Framework and applications.
\newblock In Michel Gendreau and Jean-Yves Potvin, editors, {\em Handbook of
  Metaheuristics}, volume 146 of {\em International Series in Operations
  Research \& Management Science}, pages 363--397. Springer, 2010.

\bibitem{MathWorksa}
MathWorks.
\newblock Matlab global optimization toolbox documentation.
\newblock \url{https://www.mathworks.com/help/gads/choosing-a-solver.html}.

\bibitem{MathWorks}
MathWorks.
\newblock Matlab optimization toolbox documentation.
\newblock \url{https://www.mathworks.com/help/optim/ug/choosing-a-solver.html}.

\bibitem{MerBisTra11}
Olaf Mersmann, Bernd Bischl, Heike Trautmann, Mike Preuss, Claus Weihs, and
  Günter Rudolph.
\newblock Exploratory landscape analysis.
\newblock In {\em Proceedings of the 13th annual conference on Genetic and
  evolutionary computation}, GECCO '11, pages 829--836. ACM, 2011.

\bibitem{Mersmann2013}
Olaf Mersmann, Bernd Bischl, Heike Trautmann, Markus Wagner, Jakob Bossek, and
  Frank Neumann.
\newblock A novel feature-based approach to characterize algorithm performance
  for the traveling salesperson problem.
\newblock {\em Annals of Mathematics and Artificial Intelligence},
  69(2):151--182, 2013.

\bibitem{MerPreTra11}
Olaf Mersmann, Mike Preuss, and Heike Trautmann.
\newblock Benchmarking evolutionary algorithms: Towards exploratory landscape
  analysis.
\newblock In Robert Schaefer, Carlos Cotta, Joanna Kołodziej, and Günter
  Rudolph, editors, {\em Parallel Problem Solving from Nature -- PPSN XI},
  volume 6238 of {\em Lecture Notes in Computer Science}, pages 73--82.
  Springer, 2011.

\bibitem{Miettinen2008}
Kaisa Miettinen.
\newblock Introduction to multiobjective optimization: Noninteractive
  approaches.
\newblock In Jürgen Branke, Kalyanmoy Deb, Kaisa Miettinen, and Roman
  Słowiński, editors, {\em Multiobjective Optimization}, volume 5252 of {\em
  Lecture Notes in Computer Science}, pages 1--26. Springer, 2008.

\bibitem{Mittelmann}
Hans~D. Mittelmann.
\newblock Decision tree for optimization software.
\newblock \url{http://plato.la.asu.edu/guide.html}.

\bibitem{Moen2010}
Ronald~D. Moen and Clifford~L. Norman.
\newblock Circling back.
\newblock {\em Quality Progress}, 43(11):22--28, 2010.

\bibitem{Montgomery1997}
Douglas~C. Montgomery.
\newblock {\em Design and Analysis of Experiments}.
\newblock Wiley, 4th edition, 1997.

\bibitem{Moerig2010}
Marc M{\"o}rig, Sven Scholz, Tobias Tscheuschner, and Eric Berberich.
\newblock {\em Algorithm Engineering: Bridging the Gap between Algorithm Theory
  and Practice}, chapter 6. Implementation Aspects, pages 237--289.
\newblock Springer, 2010.

\bibitem{Mueller2010}
Matthias M{\"u}ller-Hannemann and Stefan Schirra.
\newblock {\em Algorithm Engineering: Bridging the Gap between Algorithm Theory
  and Practice}, chapter 1. Foundations of Algorithm Engineering, pages 1--15.
\newblock Springer, 2010.

\bibitem{Nash1990}
John~C. Nash.
\newblock {\em Compact numerical methods for computers: Linear algebra and
  function minimisation}.
\newblock Adam Hilger Ltd., Bristol, second edition, 1990.

\bibitem{Naujoks2007}
Boris Naujoks, Max Steden, Sven-Brian Müller, and Jochen Hundemer.
\newblock Evolutionary optimization of ship propulsion systems.
\newblock In {\em IEEE Congress on Evolutionary Computation (CEC 2007)}, pages
  2809--2816, 2007.

\bibitem{Nelder1965}
John~A. Nelder and Roger Mead.
\newblock A simplex method for function minimization.
\newblock {\em The Computer Journal}, 7(4):308--313, 1965.

\bibitem{Oram2010}
Andy Oram and Greg Wilson, editors.
\newblock {\em Making Software}.
\newblock O'Reilly, 2010.

\bibitem{Preuss2007}
Mike Preuss.
\newblock Reporting on experiments in evolutionary computation.
\newblock Technical Report {CI-221/07}, University of Dortmund, Collaborative
  Research Center 531, 2007.

\bibitem{Preuss2015}
Mike Preuss, Simon Wessing, Günter Rudolph, and Gabriele Sadowski.
\newblock Solving phase equilibrium problems by means of avoidance-based
  multiobjectivization.
\newblock In Janusz Kacprzyk and Witold Pedrycz, editors, {\em Springer
  Handbook of Computational Intelligence}, pages 1159--1171. Springer, 2015.

\bibitem{Pronovost2006}
Peter Pronovost, Dale Needham, Sean Berenholtz, David Sinopoli, Haitao Chu,
  Sara Cosgrove, Bryan Sexton, Robert Hyzy, Robert Welsh, Gary Roth, Joseph
  Bander, John Kepros, and Christine Goeschel.
\newblock An intervention to decrease catheter-related bloodstream infections
  in the {ICU}.
\newblock {\em New England Journal of Medicine}, 355(26):2725--2732, 2006.
\newblock PMID: 17192537.

\bibitem{Rios2013}
Luis~Miguel Rios and Nikolaos~V. Sahinidis.
\newblock Derivative-free optimization: a review of algorithms and comparison
  of software implementations.
\newblock {\em Journal of Global Optimization}, 56(3):1247--1293, 2013.

\bibitem{Rudolph2009}
G{\"u}nter Rudolph, Mike Preuss, and Jan Quadflieg.
\newblock Two-layered surrogate modeling for tuning optimization
  metaheuristics.
\newblock In {\em ENBIS/EMSE Conference ``Design and Analysis of Computer
  Experiments'', Saint-Etienne (France), July 1-3}, 2009.
\newblock
  \url{http://www.emse.fr/enbis-emse2009/pdf/articles/Rudolph_Preuss_double-kriging-enbis2009.pdf}.

\bibitem{SacWelMit89}
Jerome Sacks, William~J. Welch, Toby~J. Mitchell, and Henry~P. Wynn.
\newblock Design and analysis of computer experiments.
\newblock {\em Statistical Science}, 4(4):409--423, 1989.

\bibitem{Sanders2009}
Peter Sanders.
\newblock Algorithm engineering -- an attempt at a definition.
\newblock In Susanne Albers, Helmut Alt, and Stefan Näher, editors, {\em
  Efficient Algorithms}, volume 5760 of {\em Lecture Notes in Computer
  Science}, pages 321--340. Springer, 2009.

\bibitem{Segura2013a}
Carlos Segura, Carlos~A. Coello~Coello, Gara Miranda, and Coromoto Le\'{o}n.
\newblock Using multi-objective evolutionary algorithms for single-objective
  optimization.
\newblock {\em 4OR}, 11(3):201--228, 2013.

\bibitem{Torczon1997}
Virginia Torczon and Michael~W. Trosset.
\newblock From evolutionary operation to parallel direct search: pattern search
  algorithms for numerical optimization.
\newblock In David~W. Scott, editor, {\em Proceedings of the 29th Symposium on
  the Interface}, volume~29 of {\em Computing Science and Statistics}, pages
  396--401, 1997.

\bibitem{TorAliVii99}
Aimo Törn, Montaz~M. Ali, and Sami Viitanen.
\newblock Stochastic global optimization: Problem classes and solution
  techniques.
\newblock {\em Journal of Global Optimization}, 14(4):437--447, 1999.

\bibitem{Vanschoren2008}
Joaquin Vanschoren, Bernhard Pfahringer, and Geoffrey Holmes.
\newblock Learning from the past with experiment databases.
\newblock In Tu-Bao Ho and Zhi-Hua Zhou, editors, {\em PRICAI 2008: Trends in
  Artificial Intelligence}, volume 5351 of {\em Lecture Notes in Computer
  Science}, pages 485--496. Springer, 2008.

\bibitem{Vanschoren2014}
Joaquin Vanschoren, Jan~N. van Rijn, Bernd Bischl, and Luis Torgo.
\newblock {OpenML}: Networked science in machine learning.
\newblock {\em SIGKDD Explorations Newsletter}, 15(2):49--60, 2014.

\bibitem{Wagner2012}
Tobias Wagner and Simon Wessing.
\newblock On the effect of response transformations in sequential parameter
  optimization.
\newblock {\em Evolutionary Computation}, 20(2):229--248, 2012.

\bibitem{Wales1999}
David~J. Wales and Harold~A. Scheraga.
\newblock Global optimization of clusters, crystals, and biomolecules.
\newblock {\em Science}, 285(5432):1368--1372, 1999.

\bibitem{Weihe2001}
Karsten Weihe.
\newblock On the differences between “practical” and “applied”.
\newblock In Stefan Näher and Dorothea Wagner, editors, {\em Algorithm
  Engineering}, volume 1982 of {\em Lecture Notes in Computer Science}, pages
  1--10. Springer, 2001.

\bibitem{Weiser2010}
Thomas~G. Weiser, Alex~B. Haynes, Gerald Dziekan, William~R. Berry, Stuart~R.
  Lipsitz, and Atul~A. Gawande.
\newblock Effect of a 19-item surgical safety checklist during urgent
  operations in a global patient population.
\newblock {\em Annals of Surgery}, 251(5):976--980, 2010.

\bibitem{Weiser2010a}
Thomas~G. Weiser, Alex~B. Haynes, Angela Lashoher, Gerald Dziekan, Daniel~J.
  Boorman, William~R. Berry, and Atul~A. Gawande.
\newblock Perspectives in quality: designing the {WHO} surgical safety
  checklist.
\newblock {\em International Journal for Quality in Health Care},
  22(5):365--370, 2010.

\bibitem{Crepinsek2014}
Matej Črepinšek, Shih-Hsi Liu, and Marjan Mernik.
\newblock Replication and comparison of computational experiments in applied
  evolutionary computing: Common pitfalls and guidelines to avoid them.
\newblock {\em Applied Soft Computing}, 19:161--170, 2014.

\end{thebibliography}

\end{document}